\documentclass{amsart}%
\usepackage[usenames]{color}
\usepackage{amsfonts}
\usepackage{graphicx}
\usepackage{amscd}
\usepackage{amsmath}
\usepackage{amssymb}
\usepackage{graphicx,pstricks,pst-plot}
\usepackage[latin1]{inputenc}
\usepackage[colorlinks=true, allcolors=blue]{hyperref}%
\providecommand{\U}[1]{\protect\rule{.1in}{.1in}}
\theoremstyle{plain}
\newtheorem{thm}{theorem}[section]
\newtheorem{theorem}[thm]{Theorem}
\newtheorem{corollary}[thm]{Corollary}
\newtheorem{lemma}[thm]{Lemma}
\newtheorem{proposition}[thm]{Proposition}

\newtheorem{definition}[thm]{Definition}
\newtheorem{remark}[thm]{Remark}
\newtheorem{example}[thm]{Example}

\begin{document}
\title[On the Number of Gradings on Matrix Algebras]{On the Number of Gradings on Matrix Algebras}
\author[D. Diniz]{Diogo Diniz}
\address{Unidade Acadêmica de Matemática, Universidade Federal de Campina Grande,
Campina Grande, PB, 58429-970, Brazil}
\email{diogo@mat.ufcg.edu.br}
\author[D. Pellegrino]{Daniel Pellegrino}
\address{Universidade Federal da Paraíba/UFPB, CCEN, Departamento de Matemática, Cidade
Universitária, 58051-900 João Pessoa, PB, Brazil}
\email{pellegrino@pq.cnpq.br and dmpellegrino@gmail.com}
\thanks{D. Diniz was partially supported by CNPq grants No.~301704/2019-8,
No.~406401/2016-0 and No.~421129/2018-2.}
\thanks{D. Pellegrino was partially supported by CNPq grant No. ~307327/2017-5 and grant
2019/0014 Paraiba State Research Foundation (FAPESQ)}

\begin{abstract}
We determine the number of isomorphism classes of elementary gradings by a
finite group on an algebra of upper block-triangular matrices. As a
consequence we prove that, for a finite abelian group $G$, the sequence of the
numbers $E(G,m)$ of isomorphism classes of elementary $G$-gradings on the algebra $M_{m}(\mathbb{F})$ of $m\times m$ matrices with
entries in a field $\mathbb{F}$  characterizes $G$. A formula for the number
of isomorphism classes of gradings by a finite abelian group on an algebra of
upper block-triangular matrices over an algebraically closed field, with mild restrictions on its characteristic, is also provided. Finally, if $G$ is a finite
abelian group, $\mathbb{F}$ is an algebraically closed field and $N(G,m)$ is the number of isomorphism classes of $G$-gradings on $M_{m}%
(\mathbb{F})$ we prove that $N(G,m)\sim\frac{1}{\left\vert G\right\vert
!}m^{\left\vert G\right\vert -1}\sim E(G,m)$.

\end{abstract}
\maketitle

\section{Introduction}

Let $\mathbb{F}$ be a field, all vector spaces, algebras and tensor products
are considered over $\mathbb{F}$. Let $A$ be an algebra and $G$ a group. A $G
$-grading on $A$ is a decomposition
\[
A=\oplus_{g\in G}A_{g}%
\]
of $A$ as a direct sum of subspaces such that $A_{g}A_{h}\subset A_{gh}$ for
all $g,h\in G$, we say that $A$ is a $G$-graded (or simply graded) algebra. If
$A=\oplus_{g\in G}A_{g}$ and $B=\oplus_{g\in G}B_{g}$ are $G$-graded algebras,
we recall that an isomorphism of $G$-graded algebras $\varphi:A\rightarrow B$
is an isomorphism of algebras such that $\varphi\left(  A_{g}\right)  \subset
B_{g}$ for all $g\in G$. We also recall that the non-zero elements $a\in A_{g}$ are called homogeneous of degree $g$.
The problem of classifying and counting all the possible $G$-gradings (up to
isomorphisms) in an algebra $A$ has been investigated by several authors. The
gradings on matrix algebras are described in \cite{BSZ01}, \cite{BZ02} and
\cite{BZ03}; we refer to the monograph \cite[Chapter 2]{EK2013} for an account
of the classification of gradings on matrix algebras. The gradings on algebras
of upper triangular matrices and on algebras of upper block-triangular
matrices are described in \cite{VZ} and \cite{Y}, respectively, and the
isomorphisms of gradings on algebras of upper block-triangular matrices are
studied in \cite{AFD} and \cite{KochY}. The problem of counting gradings on
such algebras has been investigated in \cite{r}, \cite{DKV} and \cite{KPS}.

From now on $M_{m}(\mathbb{F})$ is the algebra of $m\times m$ matrices with
entries in $\mathbb{F}$. Let $\mathbf{m}=(m_{1},\dots, m_{s})$ be an $s$-tuple
of positive integers and let $m=m_{1}+\cdots+m_{s}$. We
denote by $UT(\mathbf{m)} $ the algebra of upper block-triangular matrices in $M_{m}(\mathbb{F})$ of the form
\begin{align*}
\left(
\begin{array}
[c]{cccc}%
A_{11} & A_{12} & \cdots & A_{1s}\\
0 & A_{22} & \cdots & A_{2s}\\
\vdots & \vdots & \ddots & \vdots\\
0 & 0 & \cdots & A_{ss}%
\end{array}
\right)
\end{align*}
where $A_{ij}$ is a block of size $m_{i}\times m_{j}$.

The algebras of upper block-triangular matrices and their $\mathbb{Z}_{2}%
$-gradings appear in a fundamental way in the classification of minimal
varieties of algebras presented in \cite{GZ}, \cite{GZ2}. We remark that the
matrix algebra $M_{m}(\mathbb{F})$ and the algebra $UT_{m}(\mathbb{F})$ of
$m\times m$ upper triangular matrices over $\mathbb{F}$ are examples of
algebras of upper block-triangular matrices. A $G$-grading on $UT(\mathbf{m)}
$ is elementary if every elementary matrix $e_{ij}$ in $UT(\mathbf{m)}$ is homogeneous, we recall that $e_{ij}$ is the matrix with $1$ in the $(i,j)$-th entry and $0$ elsewhere. In
this case it is well known that there exists an $n$-tuple
$\mathbf{g}=(g_{1},\dots,g_{n})\in G^{n} $ such that $e_{ij}$ is homogeneous
of degree $g_{i}g_{j}^{-1}$. If $D$ is a $G$-graded algebra, the algebra $A=UT(\mathbf{m}%
)\otimes D$ admits a $G$-grading $A=\oplus_{g\in G}A_{g}$ such that
\[
A_{g}=\mathrm{span}\{e_{ij}\otimes d\mid d\in D_{h},g_{i}hg_{j}^{-1}=g\}.
\]

We denote by $E(G,\mathbf{m})$ the number of isomorphism classes of elementary
$G$-gradings on $UT(\mathbf{m})$ and by $N(G,\mathbf{m})$ the number of
isomorphism classes of $G$-gradings on $UT(\mathbf{m})$. For the matrix
algebra $M_{m}(\mathbb{F})$ we use the notations $E(G,m)$ and $N(G,m)$ for the
numbers of isomorphism classes of elementary and arbitrary $G$-gradings,
respectively. We provide in Theorem \ref{main}, for an arbitrary finite group
$G$, a formula for $E(G,\mathbf{m})$. A. Valenti and M. Zaizev proved that
every grading by an arbitrary group on $UT_{n}(\mathbb{F})$ is isomorphic to
an elementary grading, see \cite{VZ}. As a consequence of Theorem \ref{main}
we conclude that the number of gradings by a finite group $G$ on
$UT_{n}(\mathbb{F})$ is $\left\vert G \right\vert ^{n-1}$.
This result was obtained in \cite[Theorem 2.3]{DKV} as a consequence of the
study of graded polynomial identities for $UT_{n}(\mathbb{F})$. For finite
abelian group $G$ we provide a formula for $N(G,\mathbf{m})$ in Theorem
\ref{ngm}.

We study the asymptotic growth of $E(G,m)$ and $N(G,m)$. We prove that
$E(G,m)\sim\frac{1}{\left\vert G\right\vert !}m^{\left\vert G\right\vert -1}$
for a finite group $G$, see Proposition \ref{poly}. The main results of this
paper provide the exact asymptotic growth of $N(G,m)$ for a finite abelian
group $G$ and show that the sequence $\left(  E(G,m)\right)  _{m=1}^{\infty}$
characterizes such groups:

\begin{theorem}
\label{t1}Let $G,H$ be finite abelian groups. If $E(G,m)=E(H,m)$ for every
natural number $m$, then $G\cong H$. Moreover, in general, the result
does not hold for non-abelian groups.
\end{theorem}

The number of isomorphism classes of elementary gradings by a finite group of
prime exponent is given in Example \ref{zpton}. As a consequence
$E(G,m)=E(H,m)$, for every natural number $m$, if $G$ and $H$ are finite
groups of prime exponent and of the same order. This remark proves the final statement of Theorem \ref{t1}, i.e., the result does
not hold without the hypothesis that the groups are abelian.

\begin{theorem}
\label{t2} If the field is algebraically closed then for all finite abelian groups $G$ , the sequences of $N(G,m)$ and
$E(G,m)$ have the same asymptotic growth. More precisely, $N(G,m)\sim
E(G,m)\sim\frac{1}{\left\vert G\right\vert !}m^{\left\vert G\right\vert -1}$.
\end{theorem}

The paper is organized as follows.
In Section 2 we prove some technical results; in
Section 3 we obtain a closed formula for $E(G,\mathbf{m})$ for arbitrary
finite groups (not necessarily abelian). In particular, we provide a very
simple formula for $E(G,m)$ if $G$ is a finite group of prime exponent and we
also obtain a formula for $E(\mathbb{Z}_{n},m)$ in terms of Euler's totient function. Our first main result (Theorem \ref{t1}) is proved in the end of Section 3. In Section 4 we provide a formula for $N(G,\mathbf{m})$ for a finite abelian group $G$ and prove the second main result (Theorem \ref{t2}).

\section{Preliminaries}

In this section we present the results on the classification of gradings on
algebras of upper block-triangular matrices that will be used in the proof of
the main results of the paper.

Let $G$ be a finite group, $\mathbf{m}=(m_{1},\dots, m_{s})$ be an $s$-tuple of positive integers and
let $I_{s}=\{1,\dots, s\}$. Denote by $\gamma(\mathbf{m},G)$ the set that
consists of the maps $a:I_{s}\times G\rightarrow\mathbb{Z}$ such that
$a(i,g)\geq0$ for every $(i,g)\in I_{s}\times G$ and $\sum_{g\in
G}a(i,g)=m_{i}$ for $i=1,\dots, s$. Given $a\in\gamma(\mathbf{m}, G)$ and
$h\in G$ the map $a\cdot h$ such that $(a\cdot h)(i,g):=a(i,gh)$ lies in
$\gamma(\mathbf{m}, G) $. Note that $(a,h)\mapsto a\cdot h$ is a right action
of $G$ on $\gamma(\mathbf{m}, G)$.

Henceforth we consider this action of $G$ on $\gamma(\mathbf{m}, G)$. Our main
goal in this section is to prove that there exists a bijection from the set of
elementary gradings on $UT(\mathbf{m})$ to the set of orbits of this action.
Next we present the result in \cite[Corollary]{AFD} using the notation here
instead of the notation in \cite{AFD}, which is in terms of rings of
endomorphisms of graded flags (see \cite[Definition 1, Proposition 1]{AFD}).

Let $n=m_{1}+\cdots+m_{s}$ and let $P_{1},\dots,P_{s}$ be the subsets of
$I_{n}$ such that every element of $P_{i}$ is strictly smaller than every element of $P_{j}$ whenever $i<j$ and
$\mid P_{i}\mid=m_{i}$ for $i=1,\dots,s$. Henceforth we denote by $S_{m_{1}%
}\times\cdots\times S_{m_{s}}$ the Young subgroup associated to the partition
\begin{equation}
I_{n}=\cup_{i}P_{i},\label{112233}%
\end{equation}
of the set $I_{n}$, in this case $S_{m_{i}}=\{\sigma\in S_{n}\mid
\sigma(j)=j,\forall j\notin P_{i}\}$.

Let $D$ be an algebra with a grading by the group $G$ and let $g\in G$. We may endow the same algebra with a $G$-grading, denoted by $^{[g^{-1}]}D^{[g]}$, such that for every $h\in G$ a non-zero element $d\in D_h$ is homogeneous of degree $g^{-1}hg$ in the $G$-grading $^{[g^{-1}]}D^{[g]}$. This notation appears in the next result.

\begin{corollary}
\cite[Corollary 1]{AFD}\label{isom} Let $\mathbf{m}=(m_{1},\dots, m_{s})$ and
$\mathbf{m^{\prime}}=(m_{1}^{\prime},\dots, m_{s^{\prime}}^{\prime})$ be
tuples of positive integers and let $\mathbf{g}=(g_{1},\dots, g_{n})$ and
$\mathbf{g^{\prime}}=(g_{1},\dots, g_{n^{\prime}}^{\prime})$ be tuples of
elements of $G $, where $n=m_{1}+\cdots m_{s}$ and $n^{\prime}=m_{1}^{\prime
}+\cdots m_{s^{\prime}}^{\prime}$. Let $B$ and $B^{\prime}$ be the algebras
$UT(\mathbf{m)}$ and $UT(\mathbf{m^{\prime})}$ with the elementary gradings
induced by $\mathbf{g}$ and $\mathbf{g^{\prime}}$, respectively. Let $D$ and
$D^{\prime}$ be $G$-graded algebras with a division grading. The $G$-graded
algebras $B\otimes D$ and $B^{\prime}\otimes D^{\prime}$ are isomorphic if and
only if $\mathbf{m}=\mathbf{m^{\prime}}$, there exists a $g\in G$ such that
$^{[g^{-1}]}D^{[g]}$ is isomorphic
to $D^{\prime}$ and there exist $h_{1},\dots, h_{n}\in\mathrm{supp}\,D$ and
$\sigma\in S_{m_{1}}\times\cdots\times S_{m_{s}}$ such that $g_{i}^{\prime
}=g_{\sigma(i)}h_{\sigma(i)}g$ for $i=1,\dots, n$.
\end{corollary}

As a consequence of the corollary above we prove the main result of this section.

\begin{proposition}
\label{elcor} Let $\mathbf{m}=(m_{1},\dots, m_{s})$ be an $s$-tuple of
positive integers and let $G$ be a group. There exists a bijection between the
set of isomorphism classes of elementary $G$-gradings on $UT(\mathbf{m)}$ and
the set of orbits of the right $G$-action on $\gamma(\mathbf{m}, G)$.
\end{proposition}

\textit{Proof.} Let $n=m_{1}+\cdots+m_{s}$ and let $P_{1},\dots,P_{s}$ be the
sets in the partition (\ref{112233}) of $I_{n}$. We associate to the tuple
$\mathbf{h}=(h_{1},\dots,h_{n})\in G^{n}$ the map $a(\mathbf{h})$ such that
\[
a(\mathbf{h})(i,g):=\mid\{j\in P_{i}\mid h_{j}=g\}\mid.
\]

Note that $a(\mathbf{h})\in\gamma(\mathbf{m}, G)$ and $\mathbf{h}\mapsto
a(\mathbf{h})$ is a surjective map. If $\mathbf{g}=(g_{1},\dots, g_{n})$ and
$\mathbf{g^{\prime}}=(g_{1}^{\prime},\dots, g_{n}^{\prime})$ are tuples in
$G^{n}$ then $a(\mathbf{g})=a(\mathbf{g^{\prime}})$ if and only if there
exists $\sigma\in S_{m_{1}}\times\cdots\times S_{m_{s}}$ such that
$g_{i}^{\prime}=g_{\sigma(i)}$ for $i=1,\dots, n$. Also, $a(\mathbf{g})\cdot
g^{-1}$ is the map associated to the tuple $(g_{1}g,\dots, g_{n}g)$. Hence
$a(\mathbf{g})$ and $a(\mathbf{g^{\prime}})$ lie in the same orbit if and only
if there exist $\sigma\in S_{m_{1}}\times\cdots\times S_{m_{s}}$ and $g\in G$
such that $g_{i}^{\prime}=g_{\sigma(i)}g$ for $i=1,\dots, n$. Corollary
\ref{isom} implies that this holds if and only if the elementary gradings
induced by $\mathbf{g}$ and $\mathbf{g}^{\prime}$ are isomorphic. Therefore
the map that associates to the isomorphism class of the elementary grading on
$UT(\mathbf{m})$ induced by $\mathbf{g}\in G^{n}$ the orbit of $a(\mathbf{g})$
is a bijection. \hfill$\Box$

We end this section with the results on the classification of gradings on
algebras of upper block-triangular matrices that will be used in the proof of
Theorem \ref{ngm} that provides a formula for $N(G, \mathbf{m})$.

\begin{theorem}
\cite{Y}\label{grad} Let $G$ be any group, let $\mathbf{m}=(m_{1},m_{2}%
,\dots,m_{s})$ be an $s$-tuple of positive integers and consider any
$G$-grading the algebra $A=UT(\mathbf{m})$ of upper block-triangular matrices
over a field $\mathbb{F}$. Suppose that either $\mathrm{char}
\,\mathbb{F}=0$ or $\mathrm{char}\,\mathbb{F}>\mathrm{dim}\,A$. Then there
exists a division $G$-grading $D$ on $M_{n}(\mathbb{F})$ and an algebra
$B=UT(\mathbf{n})$, where $\mathbf{n}=(n_{1},\dots,n_{s})$, of upper
block-triangular matrices endowed with an elementary grading, such that
$A\cong B\otimes D$.
\end{theorem}

We remark that for a matrix algebra, i.e., for $s=1$ the above result holds
for an arbitrary field $\mathbb{F}$.

\begin{remark}
It follows from the known classification results of gradings on matrix
algebras that the above result holds for matrix algebras without the
hypothesis on the characteristic of the field $\mathbb{F}$, see
\cite[Corollary 2.12]{EK2013}.
\end{remark}

As commented in \cite{Y}, it is an interesting question if Theorem \ref{grad}
holds for an arbitrary field.

\section{Counting elementary gradings and the proof of Theorem \ref{t1}}%

We begin this section determining the number of isomorphism classes of
elementary gradings by a finite group $G$ on algebras of upper block-triangular matrices. We use following
well known result.

\begin{lemma}
\label{BCF} [Burnside--Cauchy--Frobenius]\label{fr} Let $G$ be a finite group
acting on a finite set $X$. Then the number of orbits equals the average
number of fixed points:
\[
\left\vert X/G\right\vert =\frac{1}{\left\vert G\right\vert } \sum
\limits_{g\in G} Fix(g),
\]
where
\[
Fix(g)=\left\vert \left\{  x\in X:x\cdot g=x\right\}  \right\vert .
\]

\end{lemma}

The number of elements of order $m$ in $\mathbb{Z}_{m}$ is given by Euler's
totient function $\phi:\mathbb{N}\rightarrow\mathbb{N}$; we recall that
$\phi(n)=n\prod_{p\mid n}\left(  1-\frac{1}{p}\right)  $. In Proposition
\ref{zn} we determine the map, in terms of $\phi$, that associates to a
positive integer $t$ the number of elements of order $t$ in $\mathbb{Z}_{m}$.
This is a particular case of the map in the next definition.

\begin{definition}
\label{mapG} Let $G$ be a finite group. We denote by $\varphi_{G}$ the map
$\mathbb{N}\rightarrow\mathbb{N}$ that associates to $t\in\mathbb{N}$ the
number of elements of $G$ of order $t$.
\end{definition}

We are now ready to provide a formula for the number of isomorphism classes of
elementary gradings by a finite group in an algebra of upper block triangular
matrices. We recall that for a finite group $G$ its exponent, denoted by
$\mathrm{Exp}(G)$, is the least common multiple of the orders of its elements.

\begin{theorem}
\label{main}
Let $G$ be a finite group and let $\mathbf{m}=(m_{1},\dots,m_{s})$ be a tuple of positive
integers. The number $E(G,\mathbf{m})$ of isomorphism classes of elementary
$G$-gradings on the algebra of upper block-triangular matrices $UT(\mathbf{m}%
)$ is
\[
E(G,\mathbf{m})=\frac{1}{\left\vert G\right\vert }\left(  \sum_{t\mid d}%
\prod_{i=1}^{s}{\binom{\frac{m_{i}}{t}+\frac{\mid G\mid}{t}-1}{\frac{m_{i}}%
{t}}}\varphi_{G}(t)\right)  ,
\]
where $d=\gcd\left(  m_{1},\dots,m_{s},\mathrm{Exp}(G)\right)$.
\end{theorem}

\textit{Proof.} We prove that for an element $h\in G$ of order $t$
\[
\mathrm{Fix}(h)=\left\{
\begin{array}
[c]{c}%
0\mbox{ if }t\nmid m_{i}\mbox{ for some }i\\ \\ 
\prod_{i=1}^{s}{\binom{q_{i}+k-1}{q_{i}}}\mbox{ if }m_{i}=q_{i}t,
\end{array}
\right.
\]
where $k=\frac{\left \vert G \right \vert}{t}$.
The result then follows from Lemma \ref{BCF}. Let $a:I_{s}\times
G\rightarrow\mathbb{Z}$ be a map such that $a\cdot h=a$. Note that
$a(i,gh)=a(i,g)$ for every $(i,g)\in I_{s}\times G$, and it follows by induction
that
\[
a(i,gh^{n})=a(i,g),
\]
for every non-negative integer $n$. Let $H=\{e,h,\dots,h^{t-1}\}$ and let
$g_{1}H,\dots,g_{u}H$ be the lateral classes of $H$ in $G$. We have
\[
m_{i}=\sum_{g\in G}a(i,g)=\sum_{j=1}^{u}\sum_{h^{\prime}\in H}a(i,g_{j}%
h^{\prime})=t\left(  \sum_{j=1}^{u}a(i,g_{j})\right)  ,
\]
therefore $t\mid m_{i}$ for $i=1,\dots,s$. Hence $\mathrm{Fix}(h)=0$ if
$t\nmid m_{i}$ for some $i\in I_{s}$. Now assume that there exist integers
$q_{1},\dots,q_{s}$ such that $m_{i}=q_{i}t$ for $i\in I_{s}$. In this case
the number of elements of $\gamma(\mathbf{m},G)$ fixed by $h$ is equal to the
number of functions $\hat{a}:I_{s}\times G/H\rightarrow\mathbb{Z}$ such that
$\hat{a}(i,gH)\geq0$ for every $(i,gH)\in I_{s}\times G/H$ and $\sum_{gH\in
G/H}\hat{a}(i,gH)=q_{i}$ for every $i\in I_{s}$. The number of such maps is
$\prod_{i=1}^{s}{\binom{q_{i}+k-1}{q_{i}}}$.

\hfill$\Box$

\vspace{0,3cm}

A particular case of the previous theorem, for matrix algebras, is the following corollary.

\begin{corollary}
\label{matrix}
Let $G$ be a finite group. For be a positive integer $m,$
the number $E(G,m)$ of isomorphism classes of elementary $G$-gradings on
$M_{m}(\mathbb{F})$ is
\[
E(G,m)=\frac{1}{\left\vert G\right\vert }\left(  \sum_{t\mid d}{\binom
{\frac{m}{t}+\frac{\mid G\mid}{t}-1}{\frac{m}{t}}}\varphi_{G}(t)\right),
\]
where $d=\gcd\left(m,\mathrm{Exp}(G)\right)$.

\end{corollary}

Next we apply of Theorem \ref{main} for the algebra of upper triangular matrices.

\begin{corollary}
	The number of isomorphism classes of gradings by a finite group $G$ on the algebra $UT_n(\mathbb{F})$ of upper triangular matrices with entries in the field $\mathbb{F}$ is $\left \vert G \right \vert^{n-1}$.
\end{corollary}
\textit{Proof.}
As a consequence of \cite[Theorem 7]{VZ} every $G$-grading on $UT_n(\mathbb{F})$ is isomorphic to an elementary grading. Hence the number of isomorphism classes of $G$-gradings on $UT_n(\mathbb{F})$ coincides with the number of isomorphism classes of elementary gradings on this algebra. The result now follows directly from Theorem \ref{main}.

\hfill $\Box$

The formula above is given in terms of the map $\varphi_{G}$ in Definition
\ref{mapG}.
If $\mathrm{Exp}(G)$ is a prime number $p$ it is plain that
\[
\varphi_{G}(t)=\left\{
\begin{array}
[c]{c}%
1\text{ if }t=1\\ \\
\left\vert G\right\vert -1\text{ if }t=p\\ \\
0\text{ if }t\nmid p
\end{array}
\right.
\]

As a consequence we obtain the following formula for $E(G,m)$ if $G$ is a
finite group of prime exponent.

\begin{example}
\label{zpton} Let $G$ be a group of prime exponent $p$ and order $p^{n}$.
Then
\[
E(G,m)=\left\{
\begin{array}
[c]{c}%
\frac{1}{p^{n}}{\binom{m+p^{n}-1}{m}}\text{ if }\gcd(p,m)=1\\ \\
\frac{1}{p^{n}}\left(  {\binom{m+p^{n}-1}{m}}+{\binom{\frac{m}{p}+p^{n-1}%
-1}{\frac{m}{p}}}(p^{n}-1)\right)  \text{ if }\gcd(p,m)=p
\end{array}
\right.
\]

\end{example}

A recurrence formula for $E(\mathbb{Z}_{p}^{n},m)$ was obtained in \cite{r},
explicit formulas are given for $n=1,2$. We remark that as a consequence of
the previous example if $G$ and $H$ are groups of prime exponent $p$ and order
$p^{n}$ then $E(G,\cdot)=E(H,\cdot)$ even if $G$ and $H$ are not isomorphic.
In Theorem \ref{t1} we prove that if $G$ and $H$ are abelian groups then
the previous equality implies that $G$ and $H$ are isomorphic.

\begin{proposition}
\label{zn} Let $n$ be a positive integer. Then
\begin{align*}
\varphi_{\mathbb{Z}_{n}}(t)=\left\{
\begin{array}
[c]{c}%
\phi(t) \text{ if } t\mid n\\ \\
0 \text{ if } t\nmid n.
\end{array}
\right.
\end{align*}

\end{proposition}

\textit{Proof.} Let $t$ be a natural number. If $t\nmid n$ then it is clear
that $\varphi_{\mathbb{Z}_{n}}(t)=0$. Now assume that $t\mid n$. The elements
of $\mathbb{Z}_{n}$ of order $t$ lie in the subgroup $H$ of $\mathbb{Z}_{n}$
of order $t$, hence $\varphi_{G}(t)=\varphi_{H}(t)$. Since $H\simeq
\mathbb{Z}_{t}$ it follows that $\varphi_{H}(t)=\phi(t)$.

\hfill$\Box$

The proposition above combined with Corollary \ref{matrix} allow us to provide a
formula for $E(\mathbb{Z}_{n},m)$ in terms of Euler's totient function. For a
cylcic group of order a power of a prime, an explicit formula is given in the
next example.

\begin{example}
\label{cicpn} Let $p$ be a prime number and let $n$ be a natural number. If
$m=p^{k}m^{\prime}$, where $p\nmid m^{\prime}$ then
\begin{align*}
E(\mathbb{Z}_{p^{n}},m)=\frac{1}{p^{n}}\left({\binom {m+p^{n}-1}{m}}+ \sum_{i=1}^{k}{\binom {p^{k-i}m^{\prime}+p^{n-i}-1}{p^{k-i}m^{\prime}}}(p^{i}-p^{i-1})\right)
,
\end{align*}

\end{example}

Note that the maps $E(\mathbb{Z}_{p}^{n}, \cdot)$ and $E(\mathbb{Z}_{p^{n}},\cdot)$
are different, however the their asymptotic growth is the same. In the next proposition we
prove that the asymptotic growth of $E(G,\cdot)$ is determined by the order of
$G$.

\begin{proposition}
\label{poly} Let $G$ be a finite group. If $d\mid\mathrm{Exp}(G)$ then there
exists a polynomial $p_{d}^{G}(x)$ of degree $\left\vert G\right\vert -1$,
leading coefficient $\frac{1}{\left\vert G\right\vert !}$ and  $p_{d}^G%
(m^{\prime})\geq0$ for every $m^{\prime}\geq0$, such that
$E(G,m)=p_{d}^{G}(m)$ whenever $\gcd(\mathrm{Exp}(G),m)=d$. Moreover $p_{1}^{G}(m)\leq
E(G,m)\leq p_{\mathrm{Exp}(G)}^{G}(m)$ for every natural number $m$, in
particular $E(G,m)\sim\frac{1}{\left\vert G\right\vert !}m^{\left\vert
G\right\vert -1}$.
\end{proposition}

\textit{Proof.} Let $f_{t}^{G}(x)=\frac{1}{\left(  \frac{\left\vert
G\right\vert }{t}-1\right)  !}(\frac{x}{t}+\frac{\left\vert G\right\vert }{t}%
-1)(\frac{x}{t}+\frac{\left\vert G\right\vert }{t}-2)\cdots(\frac{x}{t}+1)$ and let $p_{d}%
^{G}(x)=\frac{1}{ \left\vert
		G\right\vert  }\left(\sum_{t\mid d}f_{t}^G(x)\varphi_{G}(t)\right)$. Note that $p_{d}^{G}(x)$ is a
polynomial of degree $\left\vert G\right\vert -1$ and leading coefficient
$\frac{1}{\left\vert G\right\vert !}$. Corollary \ref{matrix} implies that
$E(G,m)=p_{d}^{G}(m)$ whenever $\gcd(\mathrm{Exp}(G),m)=d$. Moreover since
$f_{t}^{G}(m)>0$ for every $m$ it follows that $p_{1}^G(m)\leq E(G,m)\leq
p_{\mathrm{Exp}(G)}^G(m)$ for every natural number $m$. Since $\lim
_{m\rightarrow\infty}\frac{p_{d}(m)}{\frac{1}{\left\vert G\right\vert
!}m^{\left\vert G\right\vert !-1}}=1$ for every divisor $d$ of $\mathrm{Exp}(G)$ it follows from the previous inequalities
that $E(G,m)\sim\frac{1}{\left\vert G\right\vert !}m^{\left\vert G\right\vert
-1}$.

\hfill$\Box$

An immediate consequence of the proposition above is the following result.

\begin{corollary}
\label{asymp} If $G$ and $H$ are finite groups then $E(G,m)\sim E(H,m)$ if and
only if $\left\vert G \right\vert =\left\vert H \right\vert $.
\end{corollary}

Now we prove that the map $E(G,\cdot)$ determines $\varphi_{G}$.

\begin{proposition}
\label{varphi} Let $G$, $H$ be finite groups. The equality $E(G,m)=E(H, m)$
holds for every $m$ if and only if $\varphi_{G}=\varphi_{H}$.
\end{proposition}

\textit{Proof.} If $\varphi_{G}=\varphi_{H}$ then $\left\vert G \right\vert =
\left\vert H \right\vert $ and $\mathrm{Exp}(G)=\mathrm{Exp}(H)$; thus Corollary \ref{matrix}
implies that $E(G,m)=E(H,m)$ for every $m$. Now assume that the equality
$E(G,m)=E(H,m)$ holds for every $m$. Corollary \ref{asymp} implies that
$\left\vert G\right\vert =\left\vert H \right\vert :=k$. We have $\varphi
_{G}(1)=1=\varphi_{H}(1)$. Let $t^{\prime}$ be a natural number and assume
that $\varphi_{G}(t)=\varphi_{H}(t)$ for every $t<t^{\prime}$. We prove that
$\varphi_{G}(t^{\prime})=\varphi_{H}(t^{\prime})$, hence it follows by
induction that $\varphi_{G}(t)=\varphi_{H}(t)$ for every $t\in\mathbb{N}$. If
$t^{\prime}\nmid k $ then $\varphi_{G}(t^{\prime})=0=\varphi_{H}(t^{\prime})$.
Now assume that $t^{\prime}\mid k$. Since $E(G,t^{\prime})=E(H,t^{\prime})$,
$\left\vert G\right\vert =\left\vert H \right\vert $ and $\varphi
_{G}(t)=\varphi_{H}(t)$ for every $t<t^{\prime}$ it follows from Corollary \ref{matrix} that $\varphi_{G}(t^{\prime})=\varphi_{H}(t^{\prime})$.

\hfill$\Box$

Finally, we are able to prove Theorem \ref{t1} stated in the Introduction.

\subsection{Proof of Theorem \ref{t1}}

Proposition \ref{varphi} implies that $\varphi_{G}=\varphi_{H}$. We shall
prove that this implies that $G\cong H$. Let $p$ be a prime number and let
$G(p)$ denote the set of elements of $G$ whose order is a power of $p$. We
define $H(p)$ analogously. Corollary \ref{asymp} implies that
$\left\vert G\right\vert =\left\vert H\right\vert :=n$. For every prime number
$p$ that divides $n$, $G(p)$ and $H(p)$ are non-trivial subgroups of $G$ and
$H$, respectively and
\[
G=\oplus_{p\mid n}G(p)\text{ and }H=\oplus_{p\mid n}H(p).
\]
Note that for every prime number $p$ and every natural number $s$ we have
$\varphi_{G(p)}(p^{s})=\varphi_{G}(p^{s})$ and $\varphi_{H(p)}(p^{s}%
)=\varphi_{H}(p^{s})$. Since $\varphi_{G}=\varphi_{H}$ we conclude that
$\varphi_{G(p)}=\varphi_{H(p)}$. The result follows if we prove that
$G(p)\cong H(p)$ for every prime number $p$ that divides $n$. There exist
natural numbers $t,t^{\prime},u_{1},\dots u_{t},v_{1},\dots,v_{t^{\prime}}$ satisfying $u_{1}\geq\dots\geq u_{t}>0$ and
$v_{1}\geq\dots\geq v_{t^{\prime}}>0$ such that
\[
G(p)\cong\mathbb{Z}_{p^{u_{1}}}\times\cdots\times\mathbb{Z}_{p^{u_{t}}}\text{
and }H(p)\cong\mathbb{Z}_{p^{v_{1}}}\times\cdots\times\mathbb{Z}%
_{p^{v_{t^{\prime}}}}.
\]
Let $\psi_{G(p)}(t)$ be the number of elements of $G(p)$ of order at most $t$.
We define $\psi_{H(p)}$ analogously. Since $\varphi_{G(p)}=\varphi_{H(p)}$ we
conclude that $\psi_{G(p)}=\psi_{H(p)}$. In particular $G(p)$ and $H(p)$ have
the same order. Let $p^{\alpha}=\left\vert G(p)\right\vert $. Let $s$ be a natural
number with $1\leq s\leq u_{1}$ and let $l$ be the greatest integer such that
$u_{l}\geq s$, we have $\psi_{G(p)}(p^{s})=p^{\alpha-(u_{1}+\cdots+u_{l})+ls}$.
Note that $u_{1}$ is the greatest integer that satisfies $\psi_{G(p)}%
(p^{s})>\psi_{G}(p^{s-1})$. Since $\psi_{G(p)}=\psi_{H(p)}$ we conclude that
$u_{1}=v_{1}$. Let $r$ and $r^{\prime}$ be the greatest indexes such that
$u_{r}=u_{1}$ and $v_{r^{\prime}}=v_{1}$; for $s=u_{1}-1$ we have
\[
\psi_{G(p)}(p^{s})=p^{n-rv_{1}+rs}\text{ and }\psi_{H(p)}(p^{s}%
)=p^{n-r^{\prime}u_{1}+r^{\prime}s}.
\]
Since $\psi_{G(p)}=\psi_{H(p)}$ and $v_{1}=u_{1}$ we conclude that
$r=r^{\prime}$. We proceed in this way and conclude that $t=t^{\prime}$ and
$u_{i}=v_{i}$ for $i=1,\dots,t$. This implies that $G(p)\cong H(p)$.

\hfill$\Box$

Let $p$ be a prime number let $\alpha$ be a positive integer and let 
$\lambda=(u_{1},\dots, u_{t})\vdash\alpha$ be a partition of $\alpha$, here $u_1\geq\cdots\geq u_t>0$, and let $n=p^{\alpha}$. We denote
by $G_{p,\lambda}$ the group $\mathbb{Z}_{p^{u_{1}}}\times\cdots
\times\mathbb{Z}_{p^{u_{t}}}$. Let $\alpha_{\lambda}$ be the function given by
$\alpha_{\lambda}(s)=\alpha$ for $s>u_{1} $ and $\alpha_{\lambda}(s)=\alpha-(u_{1}%
+\cdots+ u_{l})+ls$ for $s\leq u_{1}$, where $l $ is the greatest integer such
that $u_{l}\geq s$. It follows from the proof of the previous result that
$\psi_{G_{p,\lambda}}(p^{s})=p^{\alpha_{\lambda}(s)}$. Therefore
$\varphi_{G_{p,\lambda}}(p^{s})=p^{\alpha_{\lambda}(s)}-p^{\alpha_{\lambda}(s-1)}%
$, for $s\geq 1$ and $\psi_{G_{p,\lambda}}(1)=1$. Now let $G$ be an abelian group of order $n$ and let $n=p_{1}^{\alpha_{1}%
}\cdots p_{k}^{\alpha_{k}}$ be the decomposition of $n$ as a product of prime
numbers. There exist partitions $\lambda_{i}\vdash\alpha_{i}$, for
$i=1,\dots,k$ such that $G\cong G_{p_{1},\lambda_{1}}\times\cdots\times
G_{p_{k},\lambda_{k}}$. It will be convenient to assume that
$\alpha_{\lambda(-1)}=-\infty$ and that $p^{\alpha_{\lambda}(-1)}=0$. Given
$t=p_{1}^{s_{1}}\cdots p_{k}^{s_{k}}$ be a divisor of $n$ have
\begin{align*}
\varphi_{G}(t)=\left(  p_{1}^{\alpha_{\lambda_{1}}(s_1)}-p_{1}^{\alpha
	_{\lambda_{1}}(s_1-1)}\right)  \cdots\left(  p_{k}^{\alpha_{\lambda_{k}}%
	(s_k)}-p_{k}^{\alpha_{\lambda_{k}}(s_k-1)}\right)  .
\end{align*}
Together with Corollary \ref{matrix} this allows us to obtain a formula for
$E(G,m)$ if $G$ is a finite abelian group.

\section{Counting gradings and the proof of Theorem \ref{t2}}
Our goal now is to determine the number of gradings by a finite abelian group
$G$ on an algebra of upper block-triangular matrices. We first state the
following theorem that gives the classification of division gradings on matrix
algebras, over an algebraically closed field, with support a finite abelian group.

\begin{theorem}
\cite[Theorem 2.15]{EK2013}\label{divgr} Let $T$ be a finite abelian group and
let $\mathbb{F}$ be an algebraically closed field. There exists a grading on
the matrix algebra $M_{n}(\mathbb{F})$ with support $T$ making $M_{n}%
(\mathbb{F})$ a graded division algebra if and only if $\mathrm{char}\,
\mathbb{F}$ does not divide $n$ and $T\cong\mathbb{Z}_{l_{1}}^{2}\times
\cdots\times\mathbb{Z}_{l_{r}}^{2}, l_{1}\cdots l_{r}=n$. The isomorphism
classes of such gradings are in one--to--one correspondence with the
non-degenerate alternating bicharacters $\beta:T\times T\rightarrow
\mathbb{F}^{\times}$. %All such gradings are in one equivalence class.
\end{theorem}

\begin{definition}
Let $T(G,k)=\{T\leq G \mid T\cong\mathbb{Z}_{l_{1}}^{2}\times\cdots
\times\mathbb{Z}_{l_{r}}^{2}, l_{1}\cdots l_{r}=k\}$. Let $D(T,k)$ denote the
number of isomorphism classes of division gradings with support $T$ on
$M_{k}(\mathbb{F}) $.
\end{definition}

As a consequence of Corollary \ref{isom} and the previous theorem we have the
following result.

\begin{theorem}
\label{ngm} Let $G$ be a finite abelian group. Let $\mathbf{m}=(m_{1}%
,\dots,m_{s})$ be an $s$-tuple of positive integers and let $\mathbb{F}$ be an
algebraically closed field. If $s>1$ we assume that $\mathrm{char}%
\,\mathbb{F}=0$ or $\mathrm{char}\,\mathbb{F}>\mathrm{dim}\,UT(\mathbf{m})$.
Given $k\neq0$ we denote by $\frac{\mathbf{m}}{k}$ the $s$-tuple whose $i$-th
entry is $\frac{m_{i}}{k}$. The number of isomorphism classes of $G$-gradings
on $UT(\mathbf{m)}$ is
\[
\sum_{k\mid m_{i},1\leq i\leq s}\left(  \sum_{T\in T(G,k)}D(T,k)E\left(
G/T,\frac{\mathbf{m}}{k}\right)  \right)  .
\]

\end{theorem}

\textit{Proof.} Let $(k,T,\beta, \mathbf{o})$ be a quadruple such that $k\mid
m_{i}$ for $i=1,\dots, s$, $T\in\mathcal{T}(G,k)$, $\beta:T\times
T\rightarrow\mathbb{F}^{\times}$ is a non--degenerate alternating bicharacter
and $\mathbf{o}$ is an orbit of the $G/T$-action on $\gamma(\frac{\mathbf{m}%
}{k}, G/T)$. Denote by $\mathcal{Q}$ the set of such quadruples. Let
$\mathbf{q}=(k,T,\beta, \mathbf{o})\in\mathcal{Q}$. Let $D$ be a division
grading on $M_{k}(\mathbb{F})$ that is a representative of the isomorphism
class that corresponds to $\beta:T\times T\rightarrow\mathbb{F}^{\times}$
under the one--to--one correspondence in Theorem \ref{divgr}. Let $\tilde{B}$
be an elementary grading in the class that corresponds to $\mathbf{o}$ under
the one--to--one correspondence in Proposition \ref{elcor}. Let $\mathbf{g}%
=(g_{1}T,\dots, g_{\frac{n}{k}}T)$, where $n=m_1+\cdots+m_s$, be a tuple of elements of $G/T$ that
induces the elementary grading on $\tilde{B}$ and let $B$ be the elementary
$G$-grading on $UT(\frac{\mathbf{m}}{k})$ induced by $(g_{1},\dots,
g_{\frac{n}{k}})$. The algebra $B\otimes D$ is isomorphic to $UT(\mathbf{m)}$
via the Kronecker product. Let $A$ be the corresponding $G$-grading on
$UT(\mathbf{m)}$ and denote by $[\mathbf{q}]$ the isomorphism class of $A$.
Corollary \ref{isom} implies that $[\mathbf{q}]$ is independent of the choices
of $D$, $\tilde{B}$ and $g_{1},\dots, g_{\frac{n}{k}}$. Hence we obtain a map
$\mathbf{q}\rightarrow[\mathbf{q}]$ from $\mathcal{Q}$ to the set of
isomorphism classes of $G$-gradings on $UT(\mathbf{m})$. Theorem \ref{grad}
implies that this map is surjective. Let $\mathbf{q}=(k,T,\beta, \mathbf{o})$
and $\mathbf{q}^{\prime}=(k^{\prime},T^{\prime},\beta^{\prime}, \mathbf{o}%
^{\prime})$ be quadruples in $\mathcal{Q}$ such that $[\mathbf{q}%
]=[\mathbf{q}^{\prime}]$. Let $B\otimes D$ and $B^{\prime}\otimes D^{\prime}$
be the algebras obtained from the construction above for $\mathbf{q}$ and
$\mathbf{q}^{\prime}$, respectively. Then $D$ and $D^{\prime}$ are the
division gradings on $M_{k}(\mathbb{F})$, $M_{k^{\prime}}(\mathbb{F})$
associated to $\beta$ and $\beta^{\prime}$, respectively and $\tilde{B}$,
$\tilde{B^{\prime}}$ are the elementary gradings associated to $\mathbf{o}$,
$\mathbf{o}^{\prime}$ respectively. The equality $[\mathbf{q]=[q^{\prime}]}$
implies that $B\otimes D\cong B^{\prime}\otimes D^{\prime}$. Since $G$ is
abelian we have $^{[g^{-1}]}D^{[g]}=D$ and Corollary \ref{isom} implies that
$B\otimes D\cong B^{\prime}\otimes D^{\prime}$ if and only if $D\cong
D^{\prime}$ and $\tilde{B}\cong\tilde{B^{\prime}}$. Since $D\cong D^{\prime}$ it follows
that $k=k^{\prime}$, moreover Theorem \ref{divgr} implies that $T=T^{\prime}$
and $\beta=\beta^{\prime}$. The isomorphism $\tilde{B}\cong\tilde{B^{\prime}}$ together with Proposition \ref{elcor} imply that $\mathbf{o}%
=\mathbf{o^{\prime}}$. Hence $\mathbf{q}=\mathbf{q^{\prime}}$ and the map is
also injective.

We obtain in this way a bijection from $\mathcal{Q}$ to the set of isomorphism
classes of $G$-gradings on $UT(\mathbf{m)}$. Let $k$ be a positive integer
such that $k\mid m_{i}$ for $i=1,\dots, s$ and let $T$ be a group in
$\mathcal{T}(G,k)$. Theorem \ref{divgr} implies that the number of
non-degenerate alternating bicharacters on $T$ is $D(T,k)$ and Proposition
\ref{elcor} implies that the number of orbits of the $G/T$-action on
$\gamma(\frac{\mathbf{m}}{k}, G/T)$ is $E(\frac{\mathbf{m}}{k},G/T)$. Hence
the number of quadruples $(k^{\prime},T^{\prime},\beta^{\prime},
\mathbf{o}^{\prime})\in\mathcal{Q}$ such that $k^{\prime}=k$ and $T^{\prime
}=T$ is $D(T,k)E\left(  G/T,\frac{\mathbf{m}}{k}\right)  $. As a consequence
the number of quadruples in $\mathcal{Q}$ and hence the number of isomorphism
classes of $G$-gradings on $UT(\mathbf{m})$ is $\sum_{k\mid m_{i}, 1\leq i
\leq s} \sum_{T\in T(G,k)}D(T,k)E\left(  G/T,\frac{\mathbf{m}}{k}\right)  $.

\hfill$\Box$

Now we finish the paper by proving Theorem \ref{t2}

\subsection{Proof of Theorem \ref{t2}}

Corollary \ref{ngm} implies that
\[
N(G,m)=\sum_{k\mid m}\sum_{T\in T(G,k)}D(T,k)E\left(  G/T,\frac{m}{k}\right)
.
\]
Proposition \ref{poly} implies that
\[
E(G,m)\leq N(G,m)\leq E(G,m)+\sum_{\overset{1<k\mid m}{T\in T(G,k)}}%
D(T,k)p^G_{\mathrm{Exp}(G/T)}\left(  \frac{m}{k}\right)  .
\]
The sum on the right side of the inequalities above is a polynomial of degree
strictly smaller than $\left\vert G\right\vert -1$. Hence the
above inequalities and Proposition \ref{poly} imply that $N(G,m)\sim\frac
{1}{\left\vert G\right\vert ! }m^{\left\vert G\right\vert -1}$.

\hfill$\Box$


\begin{thebibliography}{99}                                                                                               %


\bibitem {BSZ01}Y. A. Bahturin, S. K. Sehgal, and M. V. Zaicev, \textit{Group
gradings on associative algebras}, J. Algebra \textbf{241} (2001), no. 2, 677--698.

\bibitem {BZ02}Y. A. Bahturin and M. V. Zaicev, \textit{Group gradings on
matrix algebras}, Canad. Math. Bull. \textbf{45} (2002), no. 4, 499--508,
Dedicated to Robert V. Moody.

\bibitem {BZ03}Y. A. Bahturin and M. V. Zaicev, \textit{Graded algebras and
graded identities}, Polynomial identities and combinatorial methods
(Pantelleria, 2001), Lecture Notes in Pure and Appl. Math., vol. 235, Dekker,
New York, 2003, pp. 101--139.

\bibitem {r}C. Boboc, S. Dascalescu, Good gradings of matrix algebras by
finite abelian groups of prime index, Bull. Math. Soc. Sc. Math. Roumanie Tome
49 (97), No 1, 2006, 5--11.

\bibitem {AFD}\textrm{A. R. Borges, C. Fidelis and D. Diniz}, \textit{Graded
isomorphisms on upper block triangular matrix algebras}, Linear Algebra and
its Applications \textbf{543} (2018) 92--105.

\bibitem {VZ}A. Valenti, M. Zaicev, \textit{Group gradings on upper triangular
matrices}, Arch. Math. \textbf{89} (2007) 33--40.

\bibitem {DKV}Di Vincenzo, O. M., Koshlukov, P., Valenti, A., \textit{Gradings
on the algebra of upper triangular matrices and their graded identities}, J.
Algebra \textbf{275} (2004) 550--566.

\bibitem {KPS}M. Kochetov, N. Parsons, S. Sadov, Counting fine grading on
matrix algebras and on classical simple Lie algebras, Internat. J. Algebra
Comput. 23 (2013), no. 7, 1755--1781.

\bibitem {KochY}M. Kochetov, F. Yasumura, \emph{Group gradings on the Lie and
Jordan algebras of block-triangular matrices}, J. Algebra \textbf{537} (2019),
147--172.
%\bibitem {BD} M. Barascu, S. Dascalescu, Good gradings on upper block
%triangular matrix algebras, Comm. Algebra 41 (2013), no. 11, 4290--4298.


\bibitem {GZ}A. Giambruno, M. Zaicev, \emph{Codimension growth and minimal
superalgebras}, Trans. Amer. Math. Soc. \textbf{355} (2003), no. 12, 5091--5117.

\bibitem {GZ2}A. Giambruno, M. Zaicev, \emph{Minimal varieties of algebras of
exponential growth}, Adv. Math. \textbf{174} (2003), no. 2, 310--323.

\bibitem {EK2013}A. Elduque, M. Kochetov, \emph{Gradings on simple Lie
algebras}, Mathematical Surveys and Monographs, 189. American Mathematical
Society, Providence, RI; Atlantic Association for Research in the Mathematical
Sciences (AARMS), Halifax, NS, 2013.
%\bibitem {GS}A. Gordienko, O. Schnabel, On weak equivalences of gradings, J.
%Algebra 501 (2018), 435--457.


\bibitem {Y}\textrm{F.Y. Yasumura}, \textit{Group gradings on upper block
triangular matrices}, Arch. Math. \textbf{4} (2018) 327--332.
\end{thebibliography}
\end{document}